\documentclass[9pt]{article}
\usepackage{spconf,amsmath,graphicx}

\usepackage{epstopdf}

\usepackage{amsfonts,amssymb,enumerate,amsthm}
\usepackage{bbm}
\usepackage{latexsym}
\usepackage{multirow}
\usepackage{rotating}
\usepackage{mathtools}

\usepackage{algorithm}
\usepackage{algorithmic}
\usepackage{color}

\usepackage{balance}

\newcommand{\argmin}{\operatornamewithlimits{argmin}}

\title{SPARSE TENSOR RECOVERY VIA N-MODE FISTA WITH SUPPORT AUGMENTATION}
\twoauthors
  {Ashley Prater-Bennette \sthanks{This research is supported by the Air Force Office of Scientific Research under grants 18RICOR029 and 16RICOR261.  Any opinions, findings and conclusions or recommendations expressed in this material are those of the author and do not necessarily reflect the view of the United States Air Force.  Approved for public release by WPAFB on 03 July 2018, case number 88ABW-2018-3434.}}
	{Air Force Research Laboratory\\
	Information Directorate\\
	Rome, NY, USA}
  {Lixin Shen \sthanks{L.S.\ is partially supported by the US National Science Foundation under grant DMS-1522332.}}
	{Syracuse University\\
	Department of Mathematics\\
	Syracuse, NY, USA}

\begin{document}
%\ninept
%
\maketitle
\begin{abstract}
A common approach for performing sparse tensor recovery is to use an $N$-mode FISTA method.  However, this approach may fail in some cases by missing some values in the true support of the tensor and compensating by erroneously assigning nearby values to the support.  This work proposes a four-stage method for performing sparse tensor reconstruction that addresses a case where $N$-mode FISTA may fail by augmenting the support set. Moreover, the proposed method preserves a Tucker-like structure throughout computations for computational efficiency.  Numerical results on synthetic data demonstrate that the proposed method produces results with similar or higher accuracy than $N$-mode FISTA, and is often faster. 
\end{abstract}


%%%%%%%%%%%%%%%%%%%%%%%%%%%%
\section{Introduction}
Tensors are natural structures for representing multi-indexed data, and provide mechanisms for exploring relationships among several variables simultaneously in multi modal datasets \cite{kolda}.  They are essential components in myriad applications, including image and video processing~\cite{liu}, signal processing~\cite{lim}, and have been embedded within neural networks~\cite{prater_ssci}. 
Despite their strengths in representing and interpreting multidimensional data, tensors may be challenging to use in practice due to the computational issues when dealing with high dimensional data.  One may try to extend traditional matrix-vector techniques to tensors, but often these methods do not scale well.  Instead, it is desirable to develop methods that natively exploit the  multidimensional structure of tensors.

This work considers the problem of recovering a sparse core tensor $\mathcal{X} \in \mathbb{R}^{J_1 \times J_2 \times \cdots \times J_N}$ from the noisy observations
\begin{equation}\label{eq:obs}
	\mathcal{Y} = \mathcal{X} \times_1 A_1 \times_2 A_2 \times \cdots \times_N A_N + \varepsilon,
\end{equation}
where ${\mathcal{Y}\in \mathbb{R}^{I_1\times I_2 \times \cdots \times I_N}}$ with ${I_n \leq J_n}$, each ${A_n \in \mathbb{R}^{J_n \times I_n}}$ is a known factor matrix with orthonormal columns, $\varepsilon$ is a tensor of noise, and $\times_n$ refers to the $n$-mode tensor product~\cite{bader}.  Equation~\eqref{eq:obs} is a higher dimensional analogue of the problem of recovering a sparse matrix $X$ from the observations 
\begin{equation}\label{eq:cs}
	Y = A_1^\top XA_2 + \varepsilon.
\end{equation}
Problem~\eqref{eq:cs} is referred to as kronecker compressed sensing~\cite{caiafa, duarte} or sparse matrix sketching~\cite{nowak, varshney}.  
To solve Problem~\eqref{eq:cs}, the model may be recast as a classical compressed sensing matrix-vector problem
\begin{equation}\label{eq:kron cs}
	Y_\text{vec} = P X_\text{vec} + \varepsilon,
\end{equation}
where $(\cdot)_\text{vec}$ vectorizes an array by stacking its columns and $P$ is the kronecker product $P = A_2 \otimes A_1$,  with 
\begin{equation*}
	A\otimes B := \begin{bmatrix} a_{11} B & a_{12} B & \cdots &a_{1 I} B \\ a_{21}B & a_{22}B & \cdots & a_{2 I}B \\ \vdots & \vdots & \ddots & \vdots \\ a_{J 1}B & a_{J 2}B & \cdots & a_{J I}B \end{bmatrix}
\end{equation*}
for a $J\times I$ matrix $A$. 
Equation~\eqref{eq:kron cs} may then be solved using compressed sensing techniques~\cite{duarte, nowak}.

Similarly, the tensor compressed sensing problem~\eqref{eq:obs} may also be equivalently expressed in the matrix-vector format as 
\begin{equation}\label{eq:tensor kron}
	\mathcal{Y}_\text{vec} = P \mathcal{X}_\text{vec} + \varepsilon_\text{vec}
\end{equation}
with $\displaystyle P = A_N \otimes A_{N-1} \otimes \cdots A_1$.  However, even for modest sizes of $I_n$ and $J_n$, the size of the matrix $P$ may make finding a solution to~\eqref{eq:tensor kron} infeasible, as it will require storing and processing on an array of size $\prod_{n=1}^N J_n I_n$.  Utilizing a sparse recovery method that preserves the tucker tensor representation in~\eqref{eq:obs} will reduce the complexity to ${\mathcal{O}( \sum_{n=1}^N J_n I_n)}$,  representing the total size of the factor matrices $A_n$.  

In light of the above observation, this work proposes a multi-stage iterative method to approximate a solution of Problem~\eqref{eq:obs} that preserves a tensor Tucker-like structure throughout the computations. The proposed method uses an {$N$-mode} representation of the FISTA method~\cite{beck} to approximate a solution to~\eqref{eq:obs}.  If the true solution $\mathcal{X}$ does not adhere to special structure, such as block or {$N$-way} sparsity~\cite{caiafa1, nowak}, 
then FISTA may miss some support values, especially when $I \ll J$.  
This difficulty is overcome by identifying likely support values that may have been missed, then correcting the approximation on these nodes.   The proposed method is significantly faster than using Kronecker products as in  Equation~\eqref{eq:tensor kron} for most values of $J$, and is shown in Section~\ref{sec:experiments} to improve average speed and accuracy results over FISTA without the support augmentation step.  

The following notation are used in this work.  Matrices are given by uppercase letters, e.g.\ $X$, and tensors by calligraphic script, e.g.\ $\mathcal{X}$.  The individual elements in a tensor are denoted by $\mathcal{X}({\bf j})$ for indices ${\bf j} = (j_1, j_2, \ldots, j_N)$. The $\ell_1$ and Frobenius norms of tensors are higher dimensional analogues of their matrix definitions: ${\| \mathcal{X}\|_1 = \sum_{\bf j} |\mathcal{X}({\bf j})|}$ and ${\|\mathcal{X}\|_F = \left(\sum_{\bf j} |\mathcal{X}({\bf j})|^2 \right)^{1/2}}$.  A full tensor {$N$-way} product will use the shorthand notation
\begin{equation*}
	\mathcal{X}\times_{\bf n} A_{\bf n} = \mathcal{X} \times_1 A_1 \times_2 A_2 \times \cdots \times_N A_N,
\end{equation*}
or
\begin{equation*}
	\mathcal{Y}\times_{\bf n} A_{\bf n}^\top = \mathcal{Y} \times_1 A_1^\top \times_2 A_2^\top\cdots \times_N A_N^\top,
\end{equation*}
with the index ${\bf n}$ in boldfaced font on the left hand side. 

%	Figure, illustrate support miss
\begin{figure}
\centering
\includegraphics[width=0.4\textwidth]{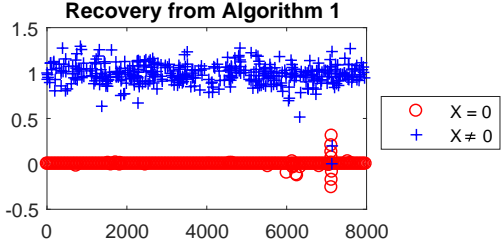}
\caption{A sample vectorized output of the $N$-mode FISTA method in Algorithm~\ref{alg:FISTA}, where a support node is assigned a value of $0$, and nearby non-support values are assigned nonzero values. Values corresponding to indices in the true support are indicated by a blue plus `${\color{blue}+}$', and values corresponding to true zero values are indicated by a red circle `${\color{red}\circ}$'.}
\label{fig:columns}
\end{figure}

%%%%%%%%%%%%%%%%%%%%%%%%%%%%%%%%%%%%%%%%%%%%%%%%%%%%%%%%%%%%%%%%%%%%%%%%%

\section{Sparse Tucker Tensor Reconstruction}
%In this section, a multi-stage method to solve Problem~\eqref{eq:obs} that preserves tensor structures throughout computations is introduced.   
%Similar to compressed sensing, a
A solution to~\eqref{eq:obs} may be approximated by a solution to the optimization problem
\begin{equation}\label{eq:tucker opt}
	\argmin_{\mathcal{U}} \left\| \mathcal{U}\right\|_1 + \frac{\lambda}{2} \left\| \mathcal{Y} - \mathcal{U}\times_{\bf n}A_{\bf n} \right\|_F^2,
\end{equation}
for some $\lambda >0$.% with $ L =  \prod_{n=1}^N \|A_n \|_2^2$.   

To solve~\eqref{eq:tucker opt}, one may use an $N$-mode iterative soft thresholding approach, such as FISTA~\cite{beck},  
or TISTA~\cite{qi}, as outlined in Algorithm~\ref{alg:FISTA}. 
Note that the $\ell_1$-proximity operator may be computed easily as a soft thresholding operator
\begin{equation}\label{eq:prox}
\mathrm{prox}_{\alpha\|\cdot\|_1}(\mathcal{U}) = \mathrm{sign}(\mathcal{U})\circ \mathrm{max}\left\{\left|\mathcal{U}\right|-\alpha,0\right\},
\end{equation}
where the operations are computed component-wise.

%	FISTA algorithm
\begin{algorithm}
	\caption{(Stage I) $N$-mode FISTA}
	\label{alg:FISTA}
	\begin{algorithmic}
	\STATE {\bfseries Inputs:} Observations $\mathcal{Y}$, factor matrices  $A_n$, parameters $\lambda, \mathrm{tol}, L>0$.\\
	\STATE {\bfseries Initialize:} Set $d_1 = 1$, $t=1$, $\mathcal{X}_0 = \mathcal{Z}_1 = \mathcal{Y}\times_{\bf n} A_{\bf n}^\top$.\\
	\STATE {\bfseries While} (stopping criteria not met) {\bfseries do} 
	\begin{align*}
		& \mathcal{Y}_t \leftarrow \mathcal{Z}_t\times_{\bf n} A_{\bf n} \\
		& \mathcal{X}_t \leftarrow \mathrm{prox}_{\frac{1}{\lambda L}\|\cdot\|_1}\left( \mathcal{Z}_t - \frac{1}{L}\left(\mathcal{Y}_t - \mathcal{Y}\right)\times_{\bf n} A_{\bf n}^\top\right) \\
		& d_{t+1} \leftarrow \frac{1+\sqrt{1+4d_t^2}}{2} \\
		&\mathcal{Z}_{t+1} \leftarrow \mathcal{X}_t + \frac{d_t-1}{d_{t+1}}\left(\mathcal{X}_{t} - \mathcal{X}_{t-1}\right) \\
		&t \leftarrow t+1
	\end{align*}
	\STATE {\bfseries End while}\\
	\STATE {\bfseries Output:} $\mathcal{X}^\text{FISTA} =$ the final $\mathcal{X}_t$ generated by the while loop.
	\end{algorithmic}
\end{algorithm}

 Although Algorithm~\ref{alg:FISTA} typically produces good results on tensors, in practice it may suffer from two shortcomings:
\begin{enumerate}
	\item Algorithm~\ref{alg:FISTA} may miss some support nodes and compensate by erroneously assigning nearby nodes to the support.  See Figure~\ref{fig:columns} for an example of this phenomenon.
	\item Common to many soft thresholding approaches, Algorithm~\ref{alg:FISTA} may underestimate values~\cite{candes, combettes, prater_dantzig}. 
\end{enumerate}

To address the first issue, the support of the approximation will be augmented as in Algorithms~\ref{alg:support} and~\ref{alg:FISTA2}.  Algorithm~\ref{alg:support} works as follows. Suppose $\widetilde{\mathcal{X}}$ is an approximate solution to Problem~\eqref{eq:tucker opt}, possibly produced by $N$-mode FISTA, with estimated support set $\widetilde{\Omega}$.  Let $\Lambda$ be a cluster of indices where $\widetilde{\mathcal{X}}$ is likely to have missed a support value.  Since any missed support values are assumed to be in a neighborhood of $\Lambda$, extend the set via
\begin{equation}\label{eq:augment}
	\widetilde{\Lambda} = \bigcup_{{\bf j}\in\Lambda} B({\bf j},r),
\end{equation}
for some radius $r>0$.  Some points in $\widetilde{\Lambda}$ may already be in the support of $\widetilde{X}$ and should not be further modified. Therefore, for some $\alpha$,  let $\mathcal{X}^\text{Aug}$ be the approximation with augmented support defined by
\begin{equation}\label{eq:augment2}
	\mathcal{X}^\text{Aug} ({\bf j}) = \begin{cases} \alpha, &\text{ if } {\bf j} \in \widetilde{\Lambda} \setminus \widetilde{\Omega}, \\ \widetilde{\mathcal{X}}({\bf j}), &\text{ else}.
	\end{cases}
\end{equation}
%for some $\alpha > b$.

%	Support augmentation
\begin{algorithm}
\caption{(Stage II) Support augmentation}
\label{alg:support}
\begin{algorithmic}
	\STATE {\bfseries Inputs:} Approximate solution $\widetilde{\mathcal{X}}$ to~\eqref{eq:tucker opt}, parameters $a, b, \gamma$, target values $\alpha$
	\STATE {\bfseries Do:}
		\begin{align*}
		&\Lambda \leftarrow \left\{ {\bf j}, {\bf k} : a < |\widetilde{\mathcal{X}}({\bf j, {\bf k}}) | < b,   \; \|{\bf j}-{\bf k}\|_F < \gamma  \right\}\\
		&\text{Compute }\widetilde{\Lambda} \text{ as in Equation~\eqref{eq:augment} } \\
		&\text{Compute } \mathcal{X}^\text{Aug} \text{ as in Equation~\eqref{eq:augment2} } \\
		& \Omega^\text{Aug} \leftarrow \widetilde{\Omega} \cup \widetilde{\Lambda}
	\end{align*}
	\STATE {\bfseries End do}
	\STATE {\bfseries Outputs:} Augmented support $\Omega^\text{Aug}$ and augmented estimate $\mathcal{X}^\text{Aug}$.
\end{algorithmic}
\end{algorithm}

Next, a modified version of FISTA is performed as in Algorithm~\ref{alg:FISTA2} that includes extra steps to project the iterative solutions onto the augmented support.  After computing $\mathcal{Y}_t$ and $\mathcal{X}_t$ as in Algorithm~\ref{alg:FISTA}, Algorithm~\ref{alg:FISTA2} updates the approximate support to include indices in the effective support of the iterative solution $\mathcal{X}_t$, then indices that have been below a certain threshold continuously for the past $R$ iterations are removed. 

%	FISTA with support
\begin{algorithm}
\caption{(Stage III) FISTA with support projections}
\label{alg:FISTA2}
\begin{algorithmic}
	\STATE {\bfseries Inputs:} Observations $\mathcal{Y}$, factor matrices $A_n$, estimate $\widetilde{\mathcal{X}}$ with support $\widetilde{\Omega}$, parameters $\lambda, \mathrm{tol}, R, L >0$.
	\STATE {\bfseries Initialize:} $ d_1 = 1, t=1, \;\mathcal{X}_0 = \mathcal{Z}_1 = \mathcal{Y}\times_{\bf n} A_{\bf n}$
	\STATE {\bfseries While} (stopping criteria not met) {\bfseries do:}
	\begin{align*}
		&\text{Compute $\mathcal{Y}_t, \mathcal{X}_t$ as in Algorithm~\ref{alg:FISTA}}\\
		&\widetilde{\Omega} \leftarrow \widetilde{\Omega} \cup \mathrm{supp}_\text{tol}(\mathcal{X}_t) \\
		& {\text {\bf if}} \left|X_s({\bf j})\right|<\mathrm{tol} \text{ for each } s = t, t-1, \ldots ,t-R+1\\
		& \quad \quad  \widetilde{\Omega} \leftarrow \widetilde{\Omega}\setminus {\bf j} \\
		&{\text {\bf end if}}\\
		&\mathcal{X}_t(\widetilde{\Omega}^c) \leftarrow {\bf 0}\\
		&\text{Update $d_{t+1}, \mathcal{Z}_{t+1},$ and $t+1$ as in Algorithm~\ref{alg:FISTA}} 
	\end{align*}
	\STATE {\bfseries End  while}
	\STATE {\bfseries Output:} $\mathcal{X}^\text{FISTA2} = $ the final $\mathcal{X}_t$ generated by the while loop.
\end{algorithmic}
\end{algorithm}

To address the second issue, a postprocessing method is used to correct the support values.  For Kronecker matrix-vector form, a least-squares postprocessing method may be used~\cite{caiafa, duarte, prater_dantzig}. Suppose $\widetilde{\mathcal{X}}$ is an approximate solution of~\eqref{eq:tucker opt} with support $\Omega$. The values of $\widetilde{\mathcal{X}}(\Omega)$ may be updated on $\Omega$ according to
\begin{equation}\label{eq:MV PP}
	\widetilde{\mathcal{X}}(\Omega) = \argmin_{x\in\mathbb{R}^{|\Omega|}} \left\{ \left\| \mathcal{Y}_\mathrm{vec} - P_\Omega x \right\|_F^2 \right\},
\end{equation}
where $P_\Omega$ is the submatrix formed by extracting the columns of $P$ with indices corresponding to $\Omega$.  As shown in Section~\ref{sec:experiments}, this approach has high accuracy, but is computationally efficient only for small sizes of $I$ and $J$.  Indeed even for modest dimensions of the factor matrices $A_n$, the kronecker product matrix $P_\Omega$ may be comptuationally infeasible to generate, store, or process with.

To this end, we propose a multi-stage method to approximation a solution to~\eqref{eq:tucker opt} that addresses the shortcomings of Algorithm~\ref{alg:FISTA}, while maintaing the tensor structure of the solution to avoid the computational issues arising from using Kronecker products.% as in Equation~\eqref{eq:MV PP}.

\begin{tabular}{p{0.52 in} p{2.5 in}}
	Stage I: 	& Run $N$-mode FISTA as in Algorithm~\ref{alg:FISTA}. \\
	Stage II:	& Augment the estimated support of the solution using Algorithm~\ref{alg:support}.\\
	Stage III:  	&Re-run $N$-mode FISTA with additional steps to project the candidate solutions onto the estimated support of $X$, as in Algorithm~\ref{alg:FISTA2}.\\
	Stage IV:  	&Post-process to refine the approximation on the augmented support set, as in Algorithm~\ref{alg:PP2}.\\
\end{tabular}

The postprocessing procedure in Stage IV, outlined in Algorithm~\ref{alg:PP2}, is similar to Algorithm~\ref{alg:FISTA}, but does not perform the soft thresholding operations and actively restricts the intermediate solutions to the estimated support set computed at the end of Stage III.  It finds an approximate solution to the tensor version of Problem~\eqref{eq:MV PP}
\begin{equation}\label{eq:tensor PP}
	\mathcal{X}^\text{PP}(\Omega) = \argmin_{\mathcal{U}} \left\| \mathcal{U} \times_{\bf n} A_{\bf n} - \mathcal{Y} \right\|_F^2 \left(\Omega\right).
\end{equation}

%	PP least squares iterations
\begin{algorithm}
\caption{(Stage IV) Iterative post processing}
\label{alg:PP2}
\begin{algorithmic}
	\STATE {\bfseries Inputs:} Observations $\mathcal{Y}$, factor matrices $A_n$, estimate $\widetilde{\mathcal{X}}$ with support $\widetilde{\Omega}$,  parameters $L, \lambda > 0$
	\STATE {\bfseries Initialize:} $d_1 = 1, t = 1, \; Z_1 = \widetilde{\mathcal{X}}$
	\STATE {\bfseries While:} (stopping criteria not met) {\bfseries do}
	\begin{align*}
		&\mathcal{Y}_t \leftarrow \mathcal{Z}_t \times_{\bf n} A_{\bf n} \\
		&\mathcal{X}_t \leftarrow \mathcal{Z}_t - \frac{\lambda }{ L} \left(\mathcal{Y}_t - \mathcal{Y}\right) \times_{\bf n} A_{\bf n}^\top  \\
		&\mathcal{X}_t (\Omega^\mathrm{c}) \leftarrow 0 \\
		&\text{Update $d_{t+1}, \mathcal{Z}_{t+1}, t+1$ as in Algorithm~\ref{alg:FISTA}}
	\end{align*}
	\STATE {\bfseries End do.}
	\STATE {\bfseries Output:}  % $\mathcal{X}^\text{PP}$, defined by
	$\displaystyle \mathcal{X}^\text{PP}({\bf j}) = \begin{cases}
		\mathcal{X}_\infty({\bf j}) 	&\text{if $|\mathcal{X}_\infty({\bf j}) | > \mathrm{tol}$ }, \\
		0 				&\text{otherwise},
		\end{cases}  $ \\
where $\mathcal{X}_\infty$ is the final $\mathcal{X}_k$ generated by the while loop. 
\end{algorithmic}
\end{algorithm}

Algorithms~\ref{alg:FISTA},~\ref{alg:FISTA2}, and~\ref{alg:PP2} may use a stopping criterion such as halting after a fixed number of iterations of once the metric $\|\mathcal{X}_t - \mathcal{X}_{t-1}\|_F$ falls below a specified threshold.

%%%%%%%%%%%%%%%%%%%%%%%%%%%%%%%%%%%%%%%%%%%%%%%%%%%%%%%%%%%%%%%%%%%%%%%%%%%%%%%

\section{Numerical Results} \label{sec:experiments}
The four-stage sparse tensor recovery method illustrated in the previous section is explored experimentally on a synthetic dataset, and is compared to three other methods: the raw $N$-mode FISTA recovery method in Algorithm~\ref{alg:FISTA}, $N$-mode FISTA with tensor postprocessing in Algorithms~\ref{alg:FISTA} and~\ref{alg:PP2}, and $N$-mode FISTA with matrix-vector postprocessing in Algorithm~\ref{alg:FISTA} and Equation~\eqref{eq:MV PP}.  The  experiments are performed in MATLAB 2017a on a PC with 2.40Gz CPU and 16GB RAM.  Notably, we do not use the Tensor Toolbox~\cite{tensor_toolbox} or CVX~\cite{cvx} to process the tensors or perform optimization, but implement the algorithms directly.

\begin{figure}[t]
\centering
\includegraphics[width=0.48\textwidth]{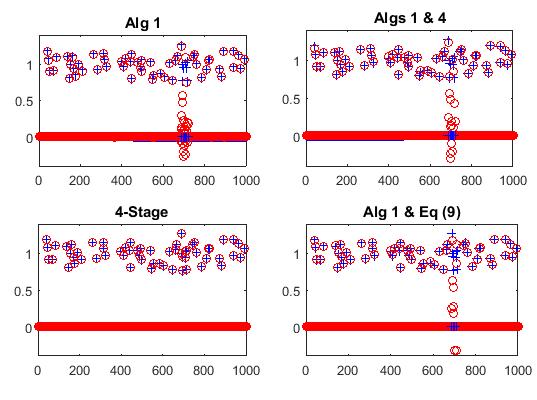}
\caption{A sample of values from one instance of a recovered tensor using each of the four methods with J = 40 and I = 28. The true values of $\mathcal{X}$ are given by a blue plus (${\color{blue} +}$), and the recovered values by a red circle (${\color{red}\circ}$).}
\label{fig:recovery}
\end{figure}

In each simulation, a sparse core tensor ${\mathcal{X}\in\mathbb{R}^{J\times J\times J}}$ is generated by randomly selecting a sparse support set $\Omega$, then sampling the entries according to $\mathcal{X}(\Omega_i) = 1 + \alpha_i,$ where each $\alpha_i$ is sampled i.i.d.\ ${\sim N(0, 0.1)}$.  The factor matrices $A_n \in\mathbb{R}^{J\times I}$, are randomly generated with orthonormal columns.  The dense tensor $\mathcal{Y}$ is observed according to
\[ \mathcal{Y} = \mathcal{X} \times_1 A_1 \times_2 A_2 \times_3 A_3 + \varepsilon,\]
where $\varepsilon$ is a noise tensor with i.i.d.\ ${\sim N(0,0.005)}$ entries.

Given the collection $\{A_n\}$ and $\mathcal{Y}$, the original sparse signal is recovered using the methods described above.  For each pair $(I,J)$ of dimensions, 20 simulations are performed with new randomizations in $\mathcal{X}$ and $A_n$ for each experiment.  However, the same randomizations are used for all recovery methods within each single experiment.  The parameters used are $\lambda = 500, \mathrm{tol} = 0.05, a = 0.05, b = 0.5, R = 20,$ and $L = 1$.

Figure~\ref{fig:recovery} displays the values from one instance of a recovered sparse tensor using the four different recovery methods.  Notice that the proposed 4-stage method in the lower left corrects the support errors that the other methods miss.

\begin{figure}
\centering
\includegraphics[width=0.48\textwidth]{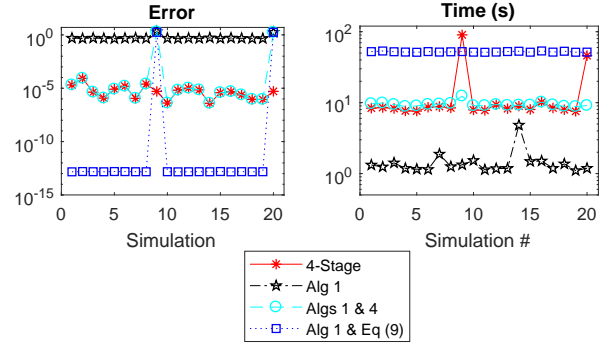}
\caption{The absolute Frobenius error (left) and time in seconds (right) to recover a sparse tensor using the four methods over 20 simulations, with $J=40$ and $I = 28$.}
\label{fig:simulation results}
\end{figure}

Figure~\ref{fig:simulation results} displays the absolute Frobenius recovery error $\|\mathcal{X} - \widetilde{\mathcal{X}}\|_F$ and the compuatation time in seconds for each recovery method computed 20 times with $I=28, J = 40$. Algorithm~\ref{alg:FISTA} did not produce the correct support in simulations 9 or 20.  In these simulations, the proposed 4-Stage method was  slower, but produced the most accurate results.  In the remaining simulations, the 4-Stage method had similar accuracy and a slight time advantage over Algorithms~\ref{alg:FISTA} and~\ref{alg:PP2}.

% Algorithm~\ref{alg:FISTA} (${\color{black}\star}$) is fast, but inaccurate.  Algorithm~\ref{alg:FISTA} postprocessed as in Equation~\ref{eq:MV PP} (${\color{blue}\square}$) is generally accurate, but slow.  The proposed 4-Stage method (${\color{red}*}$) and Algorithm~\ref{alg:FISTA} postprocessed with Algorithm~\ref{alg:PP2} (${\color{cyan}\circ}$)  in general had the same recovery error, with the 4-Stage method a bit faster.  However, on simulations when Algorithm~\ref{alg:FISTA} did not produce the correct support of the solution (simulations 9 and 20), the proposed 4-stage method was slower, but had much lower recovery error than the other methods.

\begin{figure}
\centering
\includegraphics[width=0.48\textwidth]{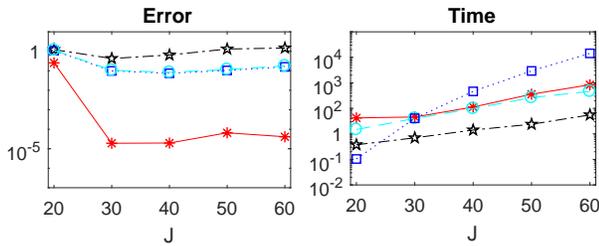}
\caption{The average error and computation time in seconds over 20 simulations to recover a sparse tensor with $I\approx 0.68J$, plotted against the values of $J$.  The same legend from Figure~\ref{fig:simulation results} is used. }
\label{fig:results}
\end{figure}

Figure~\ref{fig:results} displays the average absolute Frobenius recovery error and the computation time in seconds for several choices of $J$ with $I\approx 0.68J$. The proposed 4-Stage method achieves the lowest average recovery error, with time on average similar to that of Algorithms~\ref{alg:FISTA} and~\ref{alg:PP2}.  %When $I=J$, Algorithm~\ref{alg:FISTA} tends to produce the correct support, and therefore the 4-Stage method and Algorithm~\ref{alg:FISTA} postprocessed with Algorithm~\ref{alg:PP2} are similar in time and accuracy.    The matrix-vector postprocessing method is missing for $I=J=60$ because the Kronecker product matrix $P$ was too large to process in memory.

\section{Conclusion}

The numerical experiments illustrate that the proposed 4-Stage sparse tensor recovery method is an improvement over $N$-mode FISTA with other postprocessing approaches. The proposed method can handle instances when $N$-mode FISTA fails, while the other methods typically cannot.  Moreover, the proposed method gives comparable recovery accuracy when $N$-mode FISTA does produce the correct support. By preserving the multidimensional structure of the solution,  the proposed method reduces the computational complexity over the Kronecker matrix-vector postprocessing method.

% References should be produced using the bibtex program from suitable
% BiBTeX files (here: strings, refs, manuals). The IEEEbib.bst bibliography
% style file from IEEE produces unsorted bibliography list.
% -------------------------------------------------------------------------
\bibliographystyle{IEEEbib}
\bibliography{strings,refs}

\end{document}